\documentclass[11pt,oneside,reqno]{amsart}%
\usepackage{amssymb,amsthm,amsmath, amsfonts, float}
\usepackage{graphicx}
\usepackage{multirow}
\usepackage{subfig}
\usepackage{amsmath}
\usepackage{amsfonts}
\usepackage{lscape}
\usepackage{amssymb}%
\providecommand{\U}[1]{\protect\rule{.1in}{.1in}}
\newtheorem{theorem}{Theorem}
\theoremstyle{plain}

\newtheorem{example}{Example}

\numberwithin{equation}{section}
\ifx\pdfoutput\relax\let\pdfoutput=\undefined\fi
\newcount\msipdfoutput
\ifx\pdfoutput\undefined\else
\ifcase\pdfoutput\else
\ifx\paperwidth\undefined\else
\ifdim\paperheight=0pt\relax\else\pdfpageheight\paperheight\fi
\ifdim\paperwidth=0pt\relax\else\pdfpagewidth\paperwidth\fi
\fi\fi\fi
\begin{document}
\title[Loxodromes on Twisted Surfaces in $E_{1}^{3}$]{Loxodromes on Twisted Surfaces in Lorentz-Minkowski 3-Space}
\author[A. KAZAN]{\bfseries Ahmet KAZAN$^{1\ast}$}
\address{$^{1}$Department of Computer Technologies, Do\u{g}an\c{s}ehir Vahap
K\"{u}\c{c}\"{u}k Vocational School, Malatya Turgut \"{O}zal University,
Malatya, Turkey \\
$^{\ast}$Corresponding author: ahmet.kazan@ozal.edu.tr}
\author[M. ALTIN]{\bfseries Mustafa ALTIN$^{2}$}
\address{$^{2}$Technical Sciences Vocational School, Bing\"{o}l University, Bing\"{o}l,
Turkey }

\begin{abstract}
In the present paper, firstly we give the general formulas according to first
fundamental form of a surface for different types of loxodromes, meridians and
surfaces in $E_{1}^{3}$. After that, we obtain the differential equations of
loxodromes on Type-I, Type-II and Type-III twisted surfaces in $E_{1}^{3}$ and
also, we state a theorem which generalizes the differential equations of
different types of loxodromes on the twisted surfaces for a special case.
Finally, we provide several examples for visualizing our obtained results and
draw our loxodromes and meridians on twisted surfaces with the aid of Mathematica.

\end{abstract}
\maketitle


\section{\textbf{GENERAL INFORMATION AND BASIC CONCEPTS}}

Loxodromes (also known as rhumb lines) correspond to the curves which
intersect all of the meridians at a constant angle on the Earth. An aircraft
flying and a ship sailing on a fixed magnetic compass course move along a
curve. Here the course is a rhumb and the curve is a loxodrome. Generally, a
loxodrome is not a great circle, thus it does not measure the shortest
distance between two points on the Earth. However loxodromes are important in
navigation and they should be known by aircraft pilots and sailors \cite{Alex}.

If the shape of the Earth is approximated by a sphere, then the loxodrome is a
logarithmic spiral that cuts all meridians at the same angle and
asymptotically approaches the Earth's poles but never meets them. Since
maritime surface navigation defines the course as the angle between the
current meridian and the longitudinal direction of the ship, it may be
concluded that the loxodrome is the curve of the constant course, which means
that whenever navigating on an unchanging course we are navigating according
to a loxodrome \cite{Sergio}.

With this viewpoint, lots of studies about loxodromes in Euclidean and
Minkowskian spaces have been done in recent years. For example, the
differential equations of loxodromes on a sphere in Euclidean 3-space have
been given in \cite{Sergio} and since spheroid model of loxodrome calculation
may prove useful for understanding mathematics of navigation, the differential
equations of \ a loxodrome on the spheroid has been obtained in \cite{Miljen}.
Also, loxodromes on a twisted surface, on a canal surface, on a helicoidal
surface and on a rotational surface in $E^{3}$ have been studied in
\cite{Altin}, \cite{Baba2}, \cite{Baba1} and \cite{Kos}, respectively.
Furthermore, spacelike or timelike loxodromes on different surfaces in
Lorentz-Minkowski space have been studied in \cite{Baba4}, \cite{Baba5},
\cite{Baba3}, \cite{Baba6} and etc.

Here, firstly let we recall some basic notions about curves and twisted
surfaces in Lorentz-Minkowski 3-space $E_{1}^{3}$. For the following notions,
we refer to \cite{Dillen}, \cite{Goemans1}, \cite{oneil} and etc.

For two vectors $\overrightarrow{x}=(x_{1},$ $x_{2},$ $x_{3})$ and
$\overrightarrow{y}=(y_{1},$ $y_{2},$ $y_{3})$, let we denote the
Lorentz-Minkowski 3-space endowed with the indefinite metric
$g(\overrightarrow{x},\overrightarrow{y})=x_{1}y_{1}+x_{2}y_{2}-x_{3}y_{3}$ as
$E_{1}^{3}$. In this space, a vector $\overrightarrow{x}$ is called spacelike,
timelike and lightlike(null) if $g(\overrightarrow{x},\overrightarrow{x})$
$>0$ or $\overrightarrow{x}$ $=0$, $g(\overrightarrow{x},\overrightarrow{x})$
$<0$ and $g(\overrightarrow{x},\overrightarrow{x})$ $=0$, respectively. And
also, a curve in $E_{1}^{3}$ is said to be spacelike, timelike or
lightlike(null) when the tangent vector to the curve is everywhere spacelike,
timelike or lightlike(null), respectively. Also, the arc-length of a spacelike
or timelike curve $\alpha$, measured from $\alpha(t_{0}),$ $t_{0}\in I$ is
given by $s(t)=%
{\displaystyle\int\nolimits_{t_{0}}^{t}}
\left\Vert \dot{\alpha}(t)\right\Vert dt$. Then, the parameter $s$ is
determined as $\left\Vert \alpha^{\prime}(s)\right\Vert =1$ and here, $\alpha$
is called a unit speed curve if $\left\Vert \alpha^{\prime}(s)\right\Vert =1$
\cite{Pei}. We know that, a subspace in $E_{1}^{3}$ is timelike if it contains
a timelike vector and a subspace of $E_{1}^{3}$ is timelike if the subspace
orthogonal to it is spacelike and vice versa.

With the aid of the rotation matrices%
\[
\left(
\begin{array}
[c]{ccc}%
\cos x & -\sin x & 0\\
\sin x & \cos x & 0\\
0 & 0 & 1
\end{array}
\right)  ,\text{ \ }\left(
\begin{array}
[c]{ccc}%
1 & 0 & 0\\
0 & \cosh x & \sinh x\\
0 & \sinh x & \cosh x
\end{array}
\right)  ,\text{ \ }\left(
\begin{array}
[c]{ccc}%
\cosh x & 0 & \sinh x\\
0 & 1 & 0\\
\sinh x & 0 & \cosh x
\end{array}
\right)
\]
about the timelike $z$-axis, spacelike $x$-axis or spacelike $y$-axis,
respectively, now we can give the parametric representations of the twisted
surfaces in $E_{1}^{3}$ which are given in \cite{Goemans1}.

Firstly, let we assume that the profile curve $\alpha$ lies in the timelike
$xz$-plane. So, we can state it as $\alpha(y)=(f(y),0,g(y))$. If we rotate
this curve about the straight line through the point $(a,0,0)$ parallel with
the $y$-axis and after that rotate the obtained surface about the $z-$axis, we
can write the parametrization of twisted surface as%
\begin{equation}
\Omega(x,y)=(a+f(y)\cosh(bx)+g(y)\sinh(bx))(\cos x,\sin x,0)+(0,0,f(y)\sinh
(bx)+g(y)\cosh(bx)). \label{tip1}%
\end{equation}

Also, if we rotate the profile curve about the straight line through the point
$(0,0,a)$ parallel with the $y$-axis followed by a rotation about the
$x$-axis, we have the following parametrization of twisted surface:%
\begin{equation}
\Omega(x,y)=(a+f(y)\sinh(bx)+g(y)\cosh(bx))(0,\sinh x,\cosh x)+(f(y)\cosh
(bx)+g(y)\sinh(bx),0,0). \label{tip2}%
\end{equation}

Secondly, let we suppose that the profile curve $\alpha(y)=(f(y),g(y),0)$ lies
in the spacelike $xy$-plane. Rotating $\alpha$ about the timelike straight
line through the point $(a,0,0)$ parallel with the $z$-axis followed by a
rotation about the $y$-axis, we get the third parametrization of twisted
surface in $E_{1}^{3}$ as%
\begin{equation}
\Omega(x,y)=(a+f(y)\cos(bx)-g(y)\sin(bx))(\cosh x,0,\sinh x)+(0,f(y)\sin
(bx)+g(y)\cos(bx),0). \label{tip3}%
\end{equation}
Throuhgout the present study, we will call the above parametrizations
(\ref{tip1}), (\ref{tip2}) and (\ref{tip3}) of twisted surfaces in $E_{1}^{3}$
as Type-I, Type-II and Type-III twisted surfaces, respectively. After the
notion of twisted surface has been defined, many papers have been done in different spaces about these
surfaces in recent years. For instance in \cite{Goemans2}, twisted surfaces
with null rotation axis in $E_{1}^{3}$ have been studied and the
classification of twisted surfaces with vanishing curvature in Galilean
3-space have been given in \cite{Dede}. Also, in \cite{Kazan}, the twisted
surfaces in pseudo-Galilean space have been studied.

On the other hand, let $\overrightarrow{a}$ and $\overrightarrow{b}$ be
vectors in a Lorentzian $n$--space $E_{1}^{n}.$ Then, we have the following
cases about angles (\cite{Baba1y}, \cite{Baba2y}, \cite{Rat} and etc):

\begin{enumerate}
\item For spacelike vectors $\overrightarrow{a}$ and $\overrightarrow{b}$ that
span a spacelike vector subspace, there is a unique Lorentzian spacelike angle
$\theta\in\left[  0,\pi\right]  $ between $\overrightarrow{a}$ and
$\overrightarrow{b}$ satisfying%
\begin{equation}
\left\langle \overrightarrow{a},\overrightarrow{b}\right\rangle =\left\Vert
\overrightarrow{a}\right\Vert \left\Vert \overrightarrow{b}\right\Vert
\cos\theta. \label{a1}%
\end{equation}

\item For spacelike vectors $\overrightarrow{a}$ and $\overrightarrow{b}$ that
span a timelike vector subspace, there is a unique Lorentzian timelike angle
$\theta\in%
\mathbb{R}
^{+}$ between $\overrightarrow{a}$ and $\overrightarrow{b}$ satisfying%
\begin{equation}
\left\langle \overrightarrow{a},\overrightarrow{b}\right\rangle =\left\Vert
\overrightarrow{a}\right\Vert \left\Vert \overrightarrow{b}\right\Vert
\cosh\theta. \label{a2}%
\end{equation}

\item For a spacelike vector $\overrightarrow{a}$ and a timelike vector
$\overrightarrow{b}$ that span a timelike vector subspace, there is a unique
Lorentzian timelike angle $\theta\in%
\mathbb{R}
^{+}\cup\{0\}$ between $\overrightarrow{a}$ and $\overrightarrow{b}$
satisfying%
\begin{equation}
\left\langle \overrightarrow{a},\overrightarrow{b}\right\rangle =\left\Vert
\overrightarrow{a}\right\Vert \left\Vert \overrightarrow{b}\right\Vert
\sinh\theta. \label{a3}%
\end{equation}

\item For timelike vectors $\overrightarrow{a}$ and $\overrightarrow{b}$ that
span a timelike vector subspace, there is a unique Lorentzian timelike angle
$\theta\in%
\mathbb{R}
^{+}$ between $\overrightarrow{a}$ and $\overrightarrow{b}$ satisfying%
\begin{equation}
\left\langle \overrightarrow{a},\overrightarrow{b}\right\rangle =-\left\Vert
\overrightarrow{a}\right\Vert \left\Vert \overrightarrow{b}\right\Vert
\cosh\theta. \label{a4}%
\end{equation}

\end{enumerate}

\section{\textbf{GENERAL FORMULAS OF LOXODROMES ON A SURFACE ACCORDING TO
FIRST FUNDAMENTAL FORM IN }$E_{1}^{3}$}

Let $g_{ij},$ $i,j\in\{1,2\}$ be the coefficients of the first fundamental
form of a surface $\Omega(x,y)$ given by
\[
g_{11}=\left\langle \Omega_{x},\Omega_{x}\right\rangle ,\text{ \ }%
g_{12}=g_{21}=\left\langle \Omega_{x},\Omega_{y}\right\rangle ,\text{
\ }g_{22}=\left\langle \Omega_{y},\Omega_{y}\right\rangle .
\]
Now, we will state the general formulas of different types of loxodromes
according to the meridians and surfaces types and the surface's first
fundamental form.

A curve on a surface in $E_{1}^{3}$ which cuts all meridians ($y=$constant) or
parallels ($x=$constant) at a constant angle is called a \textit{loxodrome}.

Firstly, let we suppose that the loxodrome and meridian of a spacelike surface
are spacelike. Then, from the angle formula (\ref{a1}), we have%
\begin{equation}
\cos\theta=\frac{g_{11}dx+g_{12}dy}{\sqrt{g_{11}^{2}dx^{2}+2g_{11}%
g_{12}dxdy+g_{11}g_{22}dy^{2}}} \label{1}%
\end{equation}
and from (\ref{1}), we get%
\begin{equation}
\cos^{2}\theta\left(  g_{11}^{2}dx^{2}+2g_{11}g_{12}dxdy+g_{11}g_{22}%
dy^{2}\right)  =g_{11}^{2}dx^{2}+2g_{11}g_{12}dxdy+g_{12}^{2}dy^{2}. \label{2}%
\end{equation}
From (\ref{2}), we have%
\begin{equation}
-\sin^{2}\theta\left(  g_{11}^{2}+2g_{11}g_{12}\frac{dy}{dx}\right)  =\left(
g_{12}^{2}-g_{11}g_{22}\cos^{2}\theta\right)  \left(  \frac{dy}{dx}\right)
^{2}. \label{3}%
\end{equation}
Thus, we reach the general formula of spacelike loxodrome (with spacelike
meridian) on a spacelike surface as%
\begin{equation}
\frac{dy}{dx}=\frac{-2g_{11}g_{12}\sin^{2}\theta\mp g_{11}\sqrt{g_{11}%
g_{22}-g_{12}^{2}}\sin(2\theta)}{2\left(  g_{12}^{2}-g_{11}g_{22}\cos
^{2}\theta\right)  }. \label{4}%
\end{equation}
Now, if we assume that the loxodrome and meridian of a timelike surface are
spacelike. Then, from the angle formula (\ref{a2}), we have%
\begin{equation}
\varepsilon\cosh\theta=\frac{g_{11}dx+g_{12}dy}{\sqrt{g_{11}^{2}dx^{2}%
+2g_{11}g_{12}dxdy+g_{11}g_{22}dy^{2}}}, \label{5}%
\end{equation}
where $\varepsilon=\mp1$. Thus from (\ref{5}), we get the general formula of
spacelike loxodrome (with spacelike meridian) on a timelike surface as%
\begin{equation}
\frac{dy}{dx}=\frac{2g_{11}g_{12}\sinh^{2}\theta\mp g_{11}\sqrt{g_{12}%
^{2}-g_{11}g_{22}}\sinh(2\theta)}{2\left(  g_{12}^{2}-g_{11}g_{22}\cosh
^{2}\theta\right)  }. \label{6}%
\end{equation}
For the spacelike loxodrome and timelike meridian of a timelike surface, from
the angle formula (\ref{a3}), we get%
\begin{equation}
\varepsilon\sinh\theta=\frac{g_{11}dx+g_{12}dy}{\sqrt{-g_{11}^{2}%
dx^{2}-2g_{11}g_{12}dxdy-g_{11}g_{22}dy^{2}}}. \label{7}%
\end{equation}
Thus from (\ref{7}), we obtain the general formula of spacelike loxodrome
(with timelike meridian) on a timelike surface as%
\begin{equation}
\frac{dy}{dx}=\frac{-2g_{11}g_{12}\cosh^{2}\theta\mp g_{11}\sqrt{g_{12}%
^{2}-g_{11}g_{22}}\sinh(2\theta)}{2\left(  g_{12}^{2}+g_{11}g_{22}\sinh
^{2}\theta\right)  }. \label{8}%
\end{equation}
And lastly, if the loxodrome and meridian of a timelike surface are timelike,
from the angle formula (\ref{a4}), we have%
\begin{equation}
-\cosh\theta=\frac{g_{11}dx+g_{12}dy}{\sqrt{g_{11}^{2}dx^{2}+2g_{11}%
g_{12}dxdy+g_{11}g_{22}dy^{2}}}. \label{9}%
\end{equation}
Then from (\ref{9}), we find the general formula of timelike loxodrome (with
timelike meridian) on a timelike surface as%
\begin{equation}
\frac{dy}{dx}=\frac{2g_{11}g_{12}\sinh^{2}\theta\mp g_{11}\sqrt{g_{12}%
^{2}-g_{11}g_{22}}\sinh(2\theta)}{2\left(  g_{12}^{2}-g_{11}g_{22}\cosh
^{2}\theta\right)  }. \label{10}%
\end{equation}
Now from (\ref{4}), (\ref{6}), (\ref{8}) and (\ref{10}), we can give the
following theorem which generalizes the formulas for loxodromes on surfaces in
$E_{1}^{3}.$ Here we must note that, throughout this study we will denote,
\textquotedblleft$p=1$ or $p=-1$\textquotedblright\ (resp., \textquotedblleft%
$q=1$ or $q=-1$\textquotedblright\ and \textquotedblleft$r=1$ or
$r=-1$\textquotedblright) if the surface (resp., loxodrome and meridian) is
spacelike or timelike, respectively.

\begin{theorem}
\label{teo1}The general formulas of loxodromes on a surface in $E_{1}^{3}$ are%
\begin{equation}
\frac{dy}{dx}=\frac{-2pqrg_{11}g_{12}A^{2}\mp g_{11}\sqrt{p(g_{11}%
g_{22}-g_{12}^{2})}B}{2\left(  g_{12}^{2}-qrg_{11}g_{22}(-pA^{2}+qr)\right)
}. \label{genfor}%
\end{equation}
Here $A$ and $B$ are given by%
\begin{equation}
(A,B)=\left\{
\begin{array}
[c]{l}%
A=\sin\theta,\text{ }B=\sin(2\theta),\text{ \ \ \ \ \ if }p=1;\\
A=\sinh\theta,\text{ }B=\sinh(2\theta),\text{ \ \ if }p=-1\text{ and }qr=1;\\
A=\cosh\theta,\text{ }B=\sinh(2\theta),\text{ \ \ if }p=-1\text{ and }qr=-1.
\end{array}
\right.  \label{yyy}%
\end{equation}

\end{theorem}

\section{\textbf{LOXODROMES ON TWISTED SURFACES IN }$E_{1}^{3}$}

\subsection{Loxodromes on Type-I Twisted Surfaces in $E_{1}^{3}$}

\

Let we take the Type-I twisted surface (\ref{tip1}) parametrized by%
\begin{align}
\Omega(x,y)  &  =((a+f(y)\cosh(bx)+g(y)\sinh(bx))\cos x,\label{11}\\
&  \text{ \ \ \ \ }(a+f(y)\cosh(bx)+g(y)\sinh(bx))\sin x,f(y)\sinh
(bx)+g(y)\cosh(bx)).\nonumber
\end{align}
The coefficients of the first fundamental form of (\ref{11}) are%
\begin{equation}
\left.
\begin{array}
[c]{l}%
g_{11}=\frac{1}{2}\left(
\begin{array}
[c]{l}%
2a^{2}+\left(  1-2b^{2}+\cosh(2bx)\right)  f^{2}+\left(  2b^{2}-1+\cosh
(2bx)\right)  g^{2}\\
+4ag\sinh(bx)+4f\cosh(bx)(a+g\sinh(bx))
\end{array}
\right)  ,\\
g_{12}=g_{21}=b\left(  f^{\prime}g-fg^{\prime}\right)  ,\\
g_{22}=f^{\prime2}-g^{\prime2},
\end{array}
\right\}  \label{12}%
\end{equation}
where we denote $f=f(y),$ $g=g(y),$ $f^{\prime}=\frac{df}{dy}$ and $g^{\prime
}=\frac{dg}{dy}$.

Thus from (\ref{genfor}) and (\ref{12}), we can give the following theorem:

\begin{theorem}
The differential equations of loxodromes on the Type-I twisted surface
(\ref{11}) which cut all meridians at a constant angle $\theta$ are given by%
\begin{equation}
\frac{dy}{dx}=\frac{\left(
\begin{array}
[c]{c}%
-pqrb\left(  f^{\prime}g-fg^{\prime}\right)  A^{2}\left(
\begin{array}
[c]{l}%
2a^{2}+\left(  1-2b^{2}+\cosh(2bx)\right)  f^{2}\\
+\left(  2b^{2}-1+\cosh(2bx)\right)  g^{2}\\
+4ag\sinh(bx)+4f\cosh(bx)(a+g\sinh(bx))
\end{array}
\right) \\
\mp\frac{B}{2}\left(
\begin{array}
[c]{l}%
2a^{2}+\left(  1-2b^{2}+\cosh(2bx)\right)  f^{2}+\left(  2b^{2}-1+\cosh
(2bx)\right)  g^{2}\\
+4ag\sinh(bx)+4f\cosh(bx)(a+g\sinh(bx))
\end{array}
\right)  \times\\
\sqrt{\frac{p\left(  2b^{2}\left(  f^{\prime}g-fg^{\prime}\right)
^{2}-(f^{\prime2}-g^{\prime2})\right)  }{2}\left(
\begin{array}
[c]{l}%
2a^{2}+\left(  1-2b^{2}+\cosh(2bx)\right)  f^{2}\\
+\left(  2b^{2}-1+\cosh(2bx)\right)  g^{2}\\
+4ag\sinh(bx)+4f\cosh(bx)(a+g\sinh(bx))
\end{array}
\right)  }%
\end{array}
\right)  }{\left(
\begin{array}
[c]{c}%
2b^{2}\left(  f^{\prime}g-fg^{\prime}\right)  ^{2}\\
-qr(f^{\prime2}-g^{\prime2})(-pA^{2}+qr)
\end{array}
\right)  \left(
\begin{array}
[c]{l}%
2a^{2}+\left(  1-2b^{2}+\cosh(2bx)\right)  f^{2}\\
+\left(  2b^{2}-1+\cosh(2bx)\right)  g^{2}\\
+4ag\sinh(bx)+4f\cosh(bx)(a+g\sinh(bx))
\end{array}
\right)  }, \label{13}%
\end{equation}
where $A$ and $B$ are defined by (\ref{yyy}).
\end{theorem}

Furthermore, if we take $b=0$, from (\ref{12}) we get $g_{11}g_{22}-g_{12}%
^{2}=(a+f)^{2}\left(  f^{\prime2}-g^{\prime2}\right)  $ and so, the twisted
surface (\ref{11}) can be spacelike or timelike. Also, for the meridian
\begin{align}
\Omega(x,y_{0})  &  =(\Omega_{1}(x,y_{0}),\Omega_{2}(x,y_{0}),\Omega
_{3}(x,y_{0}))\label{14}\\
&  =((a+f(y_{0})\cosh(bx)+g(y_{0})\sinh(bx))\cos x,\nonumber\\
&  \text{ \ \ \ \ }(a+f(y_{0})\cosh(bx)+g(y_{0})\sinh(bx))\sin x,f(y_{0}%
)\sinh(bx)+g(y_{0})\cosh(bx)),\nonumber
\end{align}
we have $\left(  (\Omega_{1}(x,y_{0}))_{x}\right)  ^{2}+\left(  (\Omega
_{2}(x,y_{0}))_{x}\right)  ^{2}-\left(  (\Omega_{3}(x,y_{0}))_{x}\right)
^{2}=(a+f(y_{0}))^{2}$ and so, the meridian is always spacelike. Hence, for
the cases of our \textquotedblleft\textit{surface is spacelike, loxodrome is
spacelike, meridian is spacelike}\textquotedblright, \textquotedblleft%
\textit{surface is timelike, loxodrome is spacelike, meridian is
spacelike}\textquotedblright\ and \textquotedblleft\textit{surface is
timelike, loxodrome is timelike, meridian is spacelike}\textquotedblright\ we
can take $b=0$; but for the case of our \textquotedblleft\textit{surface is
timelike, loxodrome is spacelike, meridian is timelike}\textquotedblright\ and
\textquotedblleft\textit{surface is timelike, loxodrome is timelike, meridian
is timelike}\textquotedblright\ we cannot take $b=0.$ So, when $b=0,$ from
(\ref{4}), (\ref{6}) and (\ref{8}), the differential equation of spacelike
loxodrome (with spacelike meridian) on the spacelike Type-I twisted surface is%
\begin{equation}
\frac{\sqrt{f^{\prime2}-g^{\prime2}}}{a+f}dy=\pm\tan\theta dx, \label{15}%
\end{equation}
the differential equation of spacelike loxodrome (with spacelike meridian) on
the timelike Type-I twisted surface is%
\begin{equation}
\frac{\sqrt{-f^{\prime2}+g^{\prime2}}}{a+f}dy=\mp\tanh\theta dx \label{16}%
\end{equation}
and the differential equation of timelike loxodrome (with spacelike meridian)
on the timelike Type-I twisted surface is%
\begin{equation}
\frac{\sqrt{-f^{\prime2}+g^{\prime2}}}{a+f}dy=\pm\coth\theta dx. \label{17}%
\end{equation}

\subsection{Loxodromes on Type-II Twisted Surfaces in $E_{1}^{3}$}

\

Let we deal with the Type-II twisted surface (\ref{tip2}) parametrized by%
\begin{align}
\Omega(x,y)  &  =(f(y)\cosh(bx)+g(y)\sinh(bx),\text{\ }(a+f(y)\sinh
(bx)+g(y)\cosh(bx))\sinh x,\nonumber\\
&  \text{ \ \ \ \ }(a+f(y)\sinh(bx)+g(y)\cosh(bx))\cosh x). \label{18}%
\end{align}
The coefficients of the first fundamental form of (\ref{18}) are%
\begin{equation}
\left.
\begin{array}
[c]{l}%
g_{11}=\frac{1}{2}\left(
\begin{array}
[c]{l}%
2a^{2}+\left(  -1-2b^{2}+\cosh(2bx)\right)  f^{2}+\left(  1+2b^{2}%
+\cosh(2bx)\right)  g^{2}\\
+4ag\cosh(bx)+4f\sinh(bx)(a+g\cosh(bx))
\end{array}
\right)  ,\\
g_{12}=g_{21}=b\left(  f^{\prime}g-fg^{\prime}\right)  ,\\
g_{22}=f^{\prime2}-g^{\prime2}.
\end{array}
\right\}  \label{19}%
\end{equation}

Then from (\ref{genfor}) and (\ref{19}), we have

\begin{theorem}
The differential equations of loxodromes on the Type-II twisted surface
(\ref{18}) which cut all meridians at a constant angle $\theta$ is%
\begin{equation}
\frac{dy}{dx}=\frac{%
\begin{array}
[c]{c}%
-pqrb\left(  f^{\prime}g-fg^{\prime}\right)  A^{2}\left(
\begin{array}
[c]{l}%
2a^{2}+\left(  -1-2b^{2}+\cosh(2bx)\right)  f^{2}\\
+\left(  1+2b^{2}+\cosh(2bx)\right)  g^{2}\\
+4ag\cosh(bx)+4f\sinh(bx)(a+g\cosh(bx))
\end{array}
\right) \\
\mp\frac{B}{2}\left(
\begin{array}
[c]{l}%
2a^{2}+\left(  -1-2b^{2}+\cosh(2bx)\right)  f^{2}+\left(  1+2b^{2}%
+\cosh(2bx)\right)  g^{2}\\
+4ag\cosh(bx)+4f\sinh(bx)(a+g\cosh(bx))
\end{array}
\right)  \times\\
\sqrt{\frac{p\left(  2b^{2}\left(  f^{\prime}g-fg^{\prime}\right)
^{2}-(f^{\prime2}-g^{\prime2})\right)  }{2}\left(
\begin{array}
[c]{l}%
2a^{2}+\left(  -1-2b^{2}+\cosh(2bx)\right)  f^{2}\\
+\left(  1+2b^{2}+\cosh(2bx)\right)  g^{2}\\
+4ag\cosh(bx)+4f\sinh(bx)(a+g\cosh(bx))
\end{array}
\right)  }%
\end{array}
}{\left(
\begin{array}
[c]{c}%
2b^{2}\left(  f^{\prime}g-fg^{\prime}\right)  ^{2}\\
-qr(f^{\prime2}-g^{\prime2})(-pA^{2}+qr)
\end{array}
\right)  \left(
\begin{array}
[c]{l}%
2a^{2}+\left(  -1-2b^{2}+\cosh(2bx)\right)  f^{2}\\
+\left(  1+2b^{2}+\cosh(2bx)\right)  g^{2}\\
+4ag\cosh(bx)+4f\sinh(bx)(a+g\cosh(bx))
\end{array}
\right)  }, \label{20}%
\end{equation}
where $A$ and $B$ are defined by (\ref{yyy}).
\end{theorem}

Here, if we take $b=0$, from (\ref{19}) we get $g_{11}g_{22}-g_{12}%
^{2}=(a+g)^{2}\left(  f^{\prime2}-g^{\prime2}\right)  $ and so, our surface
can be spacelike or timelike. Also, for the meridian
\begin{align}
\Omega(x,y_{0})  &  =(\Omega_{1}(x,y_{0}),\Omega_{2}(x,y_{0}),\Omega
_{3}(x,y_{0}))\label{21}\\
&  =(f(y_{0})\cosh(bx)+g(y_{0})\sinh(bx),\text{\ }(a+f(y_{0})\sinh
(bx)+g(y_{0})\cosh(bx))\sinh x,\nonumber\\
&  \text{ \ \ \ \ }(a+f(y_{0})\sinh(bx)+g(y_{0})\cosh(bx))\cosh x),\nonumber
\end{align}
we have $\left(  (\Omega_{1}(x,y_{0}))_{x}\right)  ^{2}+\left(  (\Omega
_{2}(x,y_{0}))_{x}\right)  ^{2}-\left(  (\Omega_{3}(x,y_{0}))_{x}\right)
^{2}=(a+g(y_{0}))^{2}$ and thus, the meridian is always spacelike. Hence, for
the cases of our \textquotedblleft\textit{surface is spacelike, loxodrome is
spacelike, meridian is spacelike}\textquotedblright, \textquotedblleft%
\textit{surface is timelike, loxodrome is spacelike, meridian is
spacelike}\textquotedblright\ and \textquotedblleft\textit{surface is
timelike, loxodrome is timelike, meridian is spacelike}\textquotedblright\ we
can take $b=0$; but for the case of our \textquotedblleft\textit{surface is
timelike, loxodrome is spacelike, meridian is timelike}\textquotedblright\ and
\textquotedblleft\textit{surface is timelike, loxodrome is timelike, meridian
is timelike}\textquotedblright\ we cannot take $b=0.$ So from (\ref{4}),
(\ref{6}) and (\ref{8}), the differential equation of spacelike loxodrome
(with spacelike meridian) on the spacelike Type-II twisted surface is%
\begin{equation}
\frac{\sqrt{f^{\prime2}-g^{\prime2}}}{a+g}dy=\pm\tan\theta dx, \label{22}%
\end{equation}
the differential equation of spacelike loxodrome (with spacelike meridian) on
the timelike Type-II twisted surface is%
\begin{equation}
\frac{\sqrt{-f^{\prime2}+g^{\prime2}}}{a+g}dy=\mp\tanh\theta dx \label{23}%
\end{equation}
and the differential equation of timelike loxodrome (with spacelike meridian)
on the timelike Type-II twisted surface is%
\begin{equation}
\frac{\sqrt{-f^{\prime2}+g^{\prime2}}}{a+g}dy=\pm\coth\theta dx. \label{24}%
\end{equation}

\subsection{Loxodromes on Type-III Twisted Surfaces in $E_{1}^{3}$}

\

Let we study on the Type-III twisted surface (\ref{tip3}) parametrized by%
\begin{align}
\Omega(x,y)  &  =((a+f(y)\cos(bx)-g(y)\sin(bx))\cosh x,f(y)\sin(bx)+g(y)\cos
(bx),\nonumber\\
&  \text{ \ \ \ }(a+f(y)\cos(bx)-g(y)\sin(bx))\sinh x). \label{25}%
\end{align}
The coefficients of the first fundamental form of (\ref{25}) are%
\begin{equation}
\left.
\begin{array}
[c]{l}%
g_{11}=\frac{1}{2}\left(
\begin{array}
[c]{l}%
-2a^{2}-\left(  1-2b^{2}+\cos(2bx)\right)  f^{2}+\left(  2b^{2}-1+\cos
(2bx)\right)  g^{2}\\
+4ag\sin(bx)+4f\cos(bx)(-a+g\cosh(bx))
\end{array}
\right)  ,\\
g_{12}=g_{21}=b\left(  fg^{\prime}-f^{\prime}g\right)  ,\\
g_{22}=f^{\prime2}+g^{\prime2}.
\end{array}
\right\}  \label{26}%
\end{equation}

Thus from (\ref{genfor}) and (\ref{26}), we get

\begin{theorem}
The differential equations of loxodromes on the Type-III twisted surface
(\ref{25}) which cuts all meridians at a constant angle $\theta$ is%
\begin{equation}
\frac{dy}{dx}{\small =}\frac{%
\begin{array}
[c]{c}%
-pqrb\left(  fg^{\prime}-f^{\prime}g\right)  A^{2}\left(
\begin{array}
[c]{l}%
-2a^{2}-\left(  1-2b^{2}+\cos(2bx)\right)  f^{2}\\
+\left(  2b^{2}-1+\cos(2bx)\right)  g^{2}\\
+4ag\sin(bx)+4f\cos(bx)(-a+g\cosh(bx))
\end{array}
\right) \\
\mp\frac{B}{2}\left(
\begin{array}
[c]{l}%
-2a^{2}-\left(  1-2b^{2}+\cos(2bx)\right)  f^{2}+\left(  2b^{2}-1+\cos
(2bx)\right)  g^{2}\\
+4ag\sin(bx)+4f\cos(bx)(-a+g\cosh(bx))
\end{array}
\right)  \times\\
\sqrt{\frac{p\left(  2b^{2}\left(  fg^{\prime}-f^{\prime}g\right)
^{2}-(f^{\prime2}+g^{\prime2})\right)  }{2}\left(
\begin{array}
[c]{l}%
-2a^{2}-\left(  1-2b^{2}+\cos(2bx)\right)  f^{2}\\
+\left(  2b^{2}-1+\cos(2bx)\right)  g^{2}\\
+4ag\sin(bx)+4f\cos(bx)(-a+g\cosh(bx))
\end{array}
\right)  }%
\end{array}
}{\left(
\begin{array}
[c]{c}%
2b^{2}\left(  fg^{\prime}-f^{\prime}g\right)  ^{2}\\
-qr(f^{\prime2}+g^{\prime2})(-pA^{2}+qr)
\end{array}
\right)  \left(
\begin{array}
[c]{l}%
-2a^{2}-\left(  1-2b^{2}+\cos(2bx)\right)  f^{2}\\
+\left(  2b^{2}-1+\cos(2bx)\right)  g^{2}\\
+4ag\sin(bx)+4f\cos(bx)(-a+g\cosh(bx))
\end{array}
\right)  }, \label{27}%
\end{equation}
where $A$ and $B$ are defined by (\ref{yyy}).
\end{theorem}

Now, if we take $b=0$, from (\ref{26}) we get $g_{11}g_{22}-g_{12}%
^{2}=-(a+f)^{2}\left(  f^{\prime2}+g^{\prime2}\right)  $ and so, our surface
is always timelike. Also, for the meridian
\begin{align}
\Omega(x,y_{0})  &  =(\Omega_{1}(x,y_{0}),\Omega_{2}(x,y_{0}),\Omega
_{3}(x,y_{0}))\label{28}\\
&  =((a+f(y_{0})\cos(bx)-g(y_{0})\sin(bx))\cosh x,f(y_{0})\sin(bx)+g(y_{0}%
)\cos(bx),\nonumber\\
&  \text{ \ \ \ \ }(a+f(y_{0})\cos(bx)-g(y_{0})\sin(bx))\sinh x),\nonumber
\end{align}
we have $\left(  (\Omega_{1}(x,y_{0}))_{x}\right)  ^{2}+\left(  (\Omega
_{2}(x,y_{0}))_{x}\right)  ^{2}-\left(  (\Omega_{3}(x,y_{0}))_{x}\right)
^{2}=-(a+f(y_{0}))^{2}$ and hence, the meridian is always timelike. So, for
the cases of our \textquotedblleft\textit{surface is timelike, loxodrome is
spacelike, meridian is timelike}\textquotedblright\ and \textquotedblleft%
\textit{surface is timelike, loxodrome is timelike, meridian is timelike}%
\textquotedblright\ we can take $b=0$; but for the case of our
\textquotedblleft\textit{surface is spacelike, loxodrome is spacelike,
meridian is spacelike}\textquotedblright, \textquotedblleft\textit{surface is
timelike, loxodrome is spacelike, meridian is spacelike}\textquotedblright%
\ and \textquotedblleft\textit{surface is timelike, loxodrome is timelike,
meridian is spacelike}\textquotedblright\ we cannot take $b=0.$ So from
(\ref{8}) and (\ref{10}), the differential equation of spacelike loxodrome
(with timelike meridian) on the timelike Type-III twisted surface is%
\begin{equation}
\frac{\sqrt{f^{\prime2}+g^{\prime2}}}{a+f}dy=\mp\coth\theta dx \label{29}%
\end{equation}
and the differential equation of timelike loxodrome (with timelike meridian)
on the timelike Type-III twisted surface is%
\begin{equation}
\frac{\sqrt{f^{\prime2}+g^{\prime2}}}{a+f}dy=\pm\tanh\theta dx. \label{30}%
\end{equation}
\

Therefore, from (\ref{15})-(\ref{17}), (\ref{22})-(\ref{24}), (\ref{29}) and
(\ref{30}) we can state the following theorem which gives the differential
equations of different types of loxodromes on the Type-I, Type-II and Type-III
twisted surfaces for the case of $b=0$.

\begin{theorem}
The differential equations of loxodromes on the Type-I (expect the cases of
$p=-q=r=-1$ and $p=q=r=-1$), Type-II (expect the cases of $p=-q=r=-1$ and
$p=q=r=-1$) and Type-III (expect the cases of $p=q=r=1,$ $-p=q=r=1$ and
$-p=-q=r=1$) twisted surfaces in $E_{1}^{3}$ are%
\begin{equation}
\frac{\sqrt{p(rf^{\prime2}-g^{\prime2})}}{C}dy=\pm\frac{A}{\sqrt{qr-pA^{2}}%
}dx. \label{31}%
\end{equation}
Here, $A$ is defined by (\ref{yyy}) and%
\[
C=\left\{
\begin{array}
[c]{l}%
a+f,\text{ \ \ for Type-I and Type-III twisted surfaces;}\\
a+g,\text{ \ \ for Type-II twisted surfaces.}%
\end{array}
\right.
\]

\end{theorem}

\section{\textbf{VISUALIZATIONS FOR LOXODROMES ON TWISTED SURFACES IN }%
$E_{1}^{3}$}

In this section, we will give some examples for different types of loxodromes
on the Type-I, Type-II and Type-III twisted surfaces in $E_{1}^{3}$ and draw
their graphics with the aid of Mathematica.

\begin{example}
\label{ex1}Let we take the profile curve as $\alpha(y)=(f(y),0,g(y))=(\cosh
y,0,\sinh y)$ and $a=1,$ $b=0$ in Type-I twisted surface (\ref{tip1}). Then we
get the timelike twisted surface%
\begin{equation}
\Omega(x,y)=(\left(  1+\cosh y\right)  \cos x,\left(  1+\cosh y\right)  \sin
x,\sinh y). \label{32}%
\end{equation}
Here, our meridian is spacelike and from (\ref{6}) the differential equation
of the loxodrome is%
\begin{equation}
\frac{dy}{-(1+\cosh y)}=\tanh\theta dx. \label{34}%
\end{equation}
Solving (\ref{34}), we get $y=-2arctanh(x\tanh\theta)$ and so, our spacelike
loxodrome for $x\in(-1,1)$ and $\theta=1$ can be parametrized by%
\begin{align}
l(x)  &  =(\left(  1+\cosh(-2arctanh(x\tanh\theta))\right)  \cos
x,\label{35}\\
&  \text{ \ \ \ }\left(  1+\cosh(-2arctanh(x\tanh\theta))\right)  \sin
x,\sinh(-2arctanh(x\tanh\theta))).\nonumber
\end{align}
Also, we get $y\in(-2,2)$ and the arc-length of the loxodrome (\ref{35}) is
approximately equal to $3.40367$. In Figure 1, the meridian(red) (for $y=1$)
and the loxodrome(blue) can be seen on twisted surface (\ref{32}).
\end{example}

\begin{figure}[H]
\centering
\includegraphics[
height=1.7in, width=2.2in
]{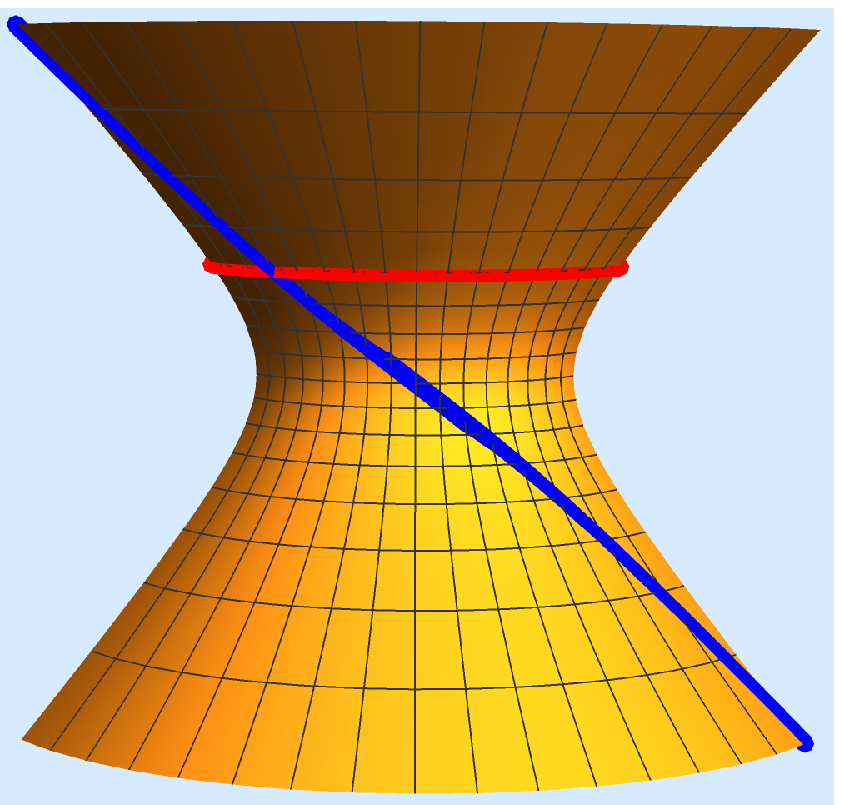}\caption{ }%
\label{fig:y}%
\end{figure}

After expressing broadly the Example \ref{ex1}, we can give the following
examples in summary form.

\begin{example}
For the profile curve $\alpha(y)=(f(y),0,g(y))=(0,0,y)$ and $a=0,$ $b=1,$ the
timelike Type-I twisted surface (\ref{tip1}) is obtained by%
\begin{equation}
\Omega(x,y)=(y\cos x\sinh x,y\sin x\sinh x,y\cosh x). \label{33}%
\end{equation}
Now, our meridian is spacelike and from (\ref{8}) the loxodrome is timelike
for $x\in(-2,2)$ and $\theta=10$. Also, we get $y\in(0.0265996,37.5946)$ and
the arc-length of the spacelike loxodrome is approximately equal to
$0.00341117$. One can see the meridian(red) (for $y=20$) and loxodrome(blue)
on twisted surface (\ref{33}) in Figure 2(a) .
\end{example}

\begin{example}
Taking $\alpha(y)=(f(y),0,g(y))=(0,0,y)$ and $a=1,$ $b=0,$ the timelike Type-I
twisted surface (\ref{tip1}) is%
\begin{equation}
\Omega(x,y)=(\cos x,\sin x,y). \label{36}%
\end{equation}
Here, our meridian is spacelike and from (\ref{8}) the loxodrome is timelike
for $x\in(-\pi,\pi)$ and $\theta=10$. Furthermore, we have $y\in
(-3.14159,3.14159)$ and the arc-length of the spacelike loxodrome is
approximately equal to $0.000570512$. The meridian(red) (for $y=1$) and
loxodrome(blue) can be seen on twisted surface (\ref{36}) in Figure 2(b).
\end{example}

\begin{figure}[H]
\centering
\includegraphics[
height=2.5in, width=4.5in
]{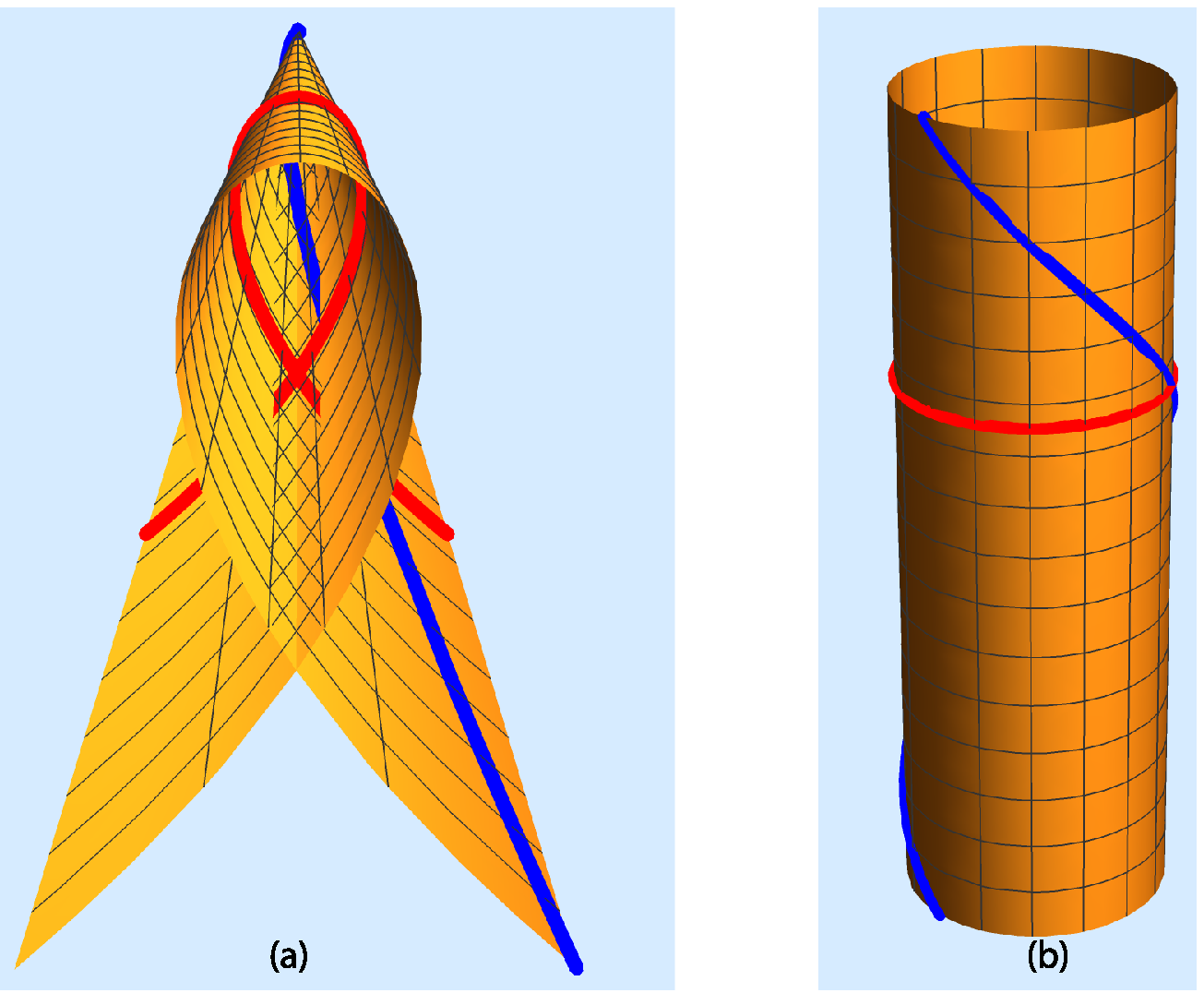}\caption{ }%
\label{fig:2}%
\end{figure}

\begin{example}
If we assume that the profile curve is $\alpha(y)=(f(y),0,g(y))=(0,0,\sinh y)$
and $a=1,$ $b=0,$ then the timelike Type-II twisted surface (\ref{tip2}) is
written by%
\begin{equation}
\Omega(x,y)=(0,(1+\sinh y)\sinh x,(1+\sinh y)\cosh x). \label{37}%
\end{equation}
Here, our meridian is spacelike and from (\ref{6}) the loxodrome is spacelike
for $x\in(-1,1)$ and $\theta=5$. Also, one can obtain that $y\in
(-0.596144,1.30993)$ and the arc-length of the spacelike loxodrome is
approximately equal to $0.0316714$. In Figure 3(a), one can see the
meridian(red) (for $y=0.8$) and loxodrome(blue) on twisted surface (\ref{37}).
\end{example}

\begin{example}
For $\alpha(y)=(f(y),0,g(y))=(y,0,0)$ and $a=0,$ $b=10,$ the timelike Type-II
twisted surface (\ref{tip2}) is%
\begin{equation}
\Omega(x,y)=(y\cosh(10x),y\sinh x\sinh(10x),y\cosh x\sinh(10x)) \label{38}%
\end{equation}
and so, our meridian is timelike and from (\ref{10}) the loxodrome is timelike
for $x\in(-0.2,0.2)$ and $\theta=10$. Here, we get $y\in(-0.139407,7.17324)$
and we obtain that the arc-length of the timelike loxodrome is approximately
equal to $0.000638671$. The meridian(red) (for $y=4$) and loxodrome(blue) can
be seen on twisted surface (\ref{38}) in Figure 3(b).
\end{example}

\begin{figure}[H]
\centering
\includegraphics[
height=1.8in, width=5.1in
]{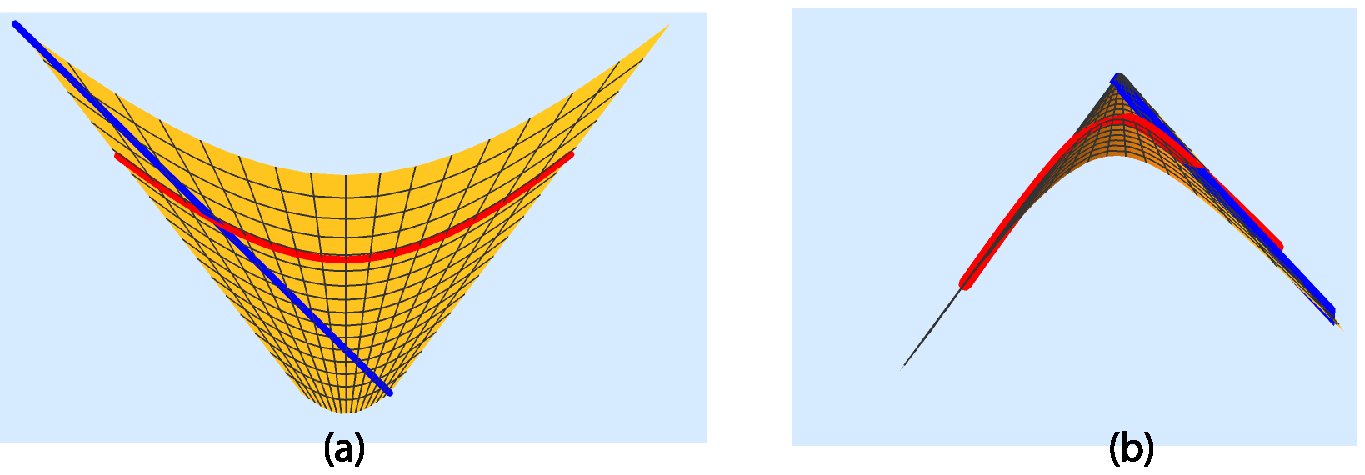}\caption{ }%
\label{fig:3}%
\end{figure}

\begin{example}
Putting $f(y)=0$ and $g(y)=y$ in the profile curve $\alpha(y)=(f(y),g(y),0)$
and taking $a=0,$ $b=1,$ the spacelike Type-III twisted surface (\ref{tip3})
is given by%
\begin{equation}
\Omega(x,y)=(-y\cosh x\sin x,y\cos x,-y\sin x\sinh x). \label{39}%
\end{equation}
Now, our meridian is spacelike and from (\ref{4}) the loxodrome is spacelike
for $x\in(-\frac{\pi}{3},\frac{\pi}{3})$ and $\theta=\frac{\pi}{4}$. Also, we
obtain $y\in(0.42062,2.37744)$ and the arc-length of the spacelike loxodrome
is approximately equal to $2.76737$. In Figure 4(a), the meridian(red) (for
$y=1$) and loxodrome(blue) can be seen on twisted surface (\ref{39}).
\end{example}

\begin{example}
Let we assume that the profile curve is $\alpha(y)=(f(y),g(y),0)=(\cosh
y,1,0)$ and $a=1,$ $b=0.$ Then, the timelike Type-III twisted surface
(\ref{tip3}) is obtained by%
\begin{equation}
\Omega(x,y)=((1+\cosh y)\cosh x,1,(1+\cosh y)\sinh x) \label{40}%
\end{equation}
and our meridian is timelike. Also from (\ref{8}), the loxodrome is spacelike
for $x\in(-1.01,-0.1)$ and $\theta=1$. Furthermore, we have $y\in
(-2.56413,-0.732672)$ and the arc-length of the spacelike loxodrome is
approximately equal to $3.40392$. In Figure 4(b), one can see the
meridian(red) (for $y=-2$) and loxodrome(blue) on twisted surface (\ref{40}).
\end{example}

\begin{figure}[H]
\centering
\includegraphics[
height=1.8in, width=5.1in
]{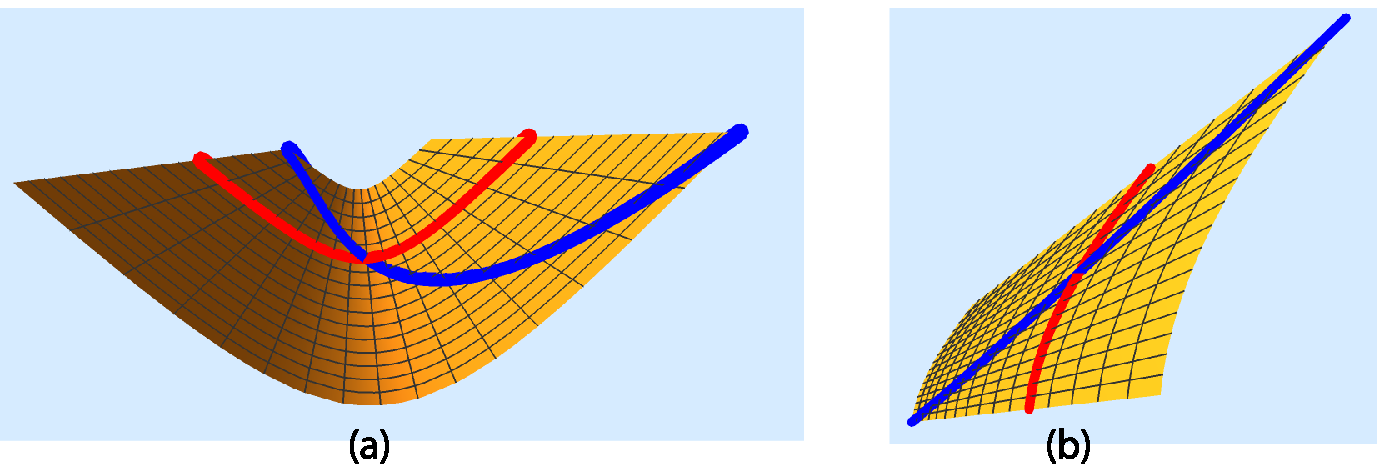}\caption{ }%
\label{fig:4}%
\end{figure}

\end{document}